\documentclass{amsart}

\input xy
\xyoption{all}
\usepackage{lacromay}
\usepackage{notation}
\usepackage{autoref_long}
\usepackage{tikz}
\usepackage{array}

\usetikzlibrary{scopes,calc,arrows,patterns}
\tikzset{ed/.style={auto,inner sep=0pt,font=\scriptsize}} 
\tikzset{>=stealth'}
\tikzset{vert/.style={draw,circle,inner sep=1pt,fill=white}}
\tikzset{vert2/.style={draw,circle,inner sep=2pt,fill=white}}

\colorlet{myblue}{blue!40!white}
\colorlet{myred}{red!35!white}
\colorlet{mygreen}{green!30!white}
\colorlet{myyellow}{yellow!10!white}
\tikzset{bluefill/.style={fill=myblue}}
\tikzset{redfill/.style={fill=myred}}
\tikzset{greenfill/.style={fill=mygreen}}
\tikzset{yellowfill/.style={fill=myyellow}}
\tikzset{dotsF/.style={pattern=north east lines,pattern color=black!60!white}}
\tikzset{dotsG/.style={pattern=north west lines,pattern color=black!60!white}}

\colorlet{mydarkred}{red!80!black}
\colorlet{mydarkblue}{blue!70!black}
\gdef\setinnerouter{\tikzset{inner/.style={densely dotted}}
  \tikzset{outer/.style={densely dashed}}
  \tikzset{outervert/.style={solid,vert}}}

\setinnerouter
\tikzset{innervert/.style={inner,circle,fill,inner sep=1pt}}
\tikzset{topvert/.style={innervert}}

\tikzset{innerunit/.style={inner,regular polygon,regular polygon sides=3,shape border rotate=180,fill,inner sep=1pt}}
\tikzset{innercounit/.style={inner,regular polygon,regular polygon sides=3,fill,inner sep=1pt}}
\tikzset{inneriso/.style={inner,circle,fill,inner sep=1pt}}
\tikzset{innerbc/.style={inner,diamond,fill,inner sep=1.3pt}}

\usetikzlibrary{decorations.pathmorphing}
\usetikzlibrary{decorations}
\tikzset{transf/.style={decorate,decoration={zigzag,amplitude=1pt,segment length=3pt}}}

\pgfdeclarelayer{background}
\pgfdeclarelayer{foreground}
\pgfsetlayers{background,main,foreground}

\def\bgcylinder#1#2#3#4#5#6{
  \def\cylempty{}\def\cylfrontcolor{#5}\def\cylbackcolor{#6}
  \begin{pgfonlayer}{background}
    \ifx\cylfrontcolor\cylempty\draw\else\fill[fill=my#5]\fi
    (#1) coordinate (dl)
    -- ++(0,#2) node[coordinate] (ul) {} 
    arc (-180:0:#3 and #4) coordinate (ur)
    -- ++(0,-#2) node[coordinate] (dr) {}
    arc (0:-180:#3 and #4);
    \ifx\cylbackcolor\cylempty\draw\else\fill[fill=my#6!80!black]\fi
    ($(ul)!.5!(ur)$) node[coordinate] (top) {} ellipse (#3 and #4);
    \path (dl) arc (-180:-90:#3 and #4) node[coordinate] (bot) {};
    \clip (dl) -- (ul) arc (-180:0:#3 and #4) -- (dr) arc (0:-180:#3 and #4);
  \end{pgfonlayer}
  \begin{pgfonlayer}{foreground}
    \clip (dl) -- (ul) arc (-180:0:#3 and #4) -- (dr) arc (0:-180:#3 and #4);
  \end{pgfonlayer}
  \clip (dl) -- (ul) arc (-180:0:#3 and #4) -- (dr) arc (0:-180:#3 and #4);
  \path (ul) ++(-0.1,0.1) coordinate (ul');
  \path (ur) ++(0.1,0.1) coordinate (ur');
  \path (dl) ++(-0.1,-0.1) coordinate (dl');
  \path (dr) ++(0.1,-0.1) coordinate (dr');
  \path (top) ++(0,0.1) coordinate (top');
  \path (bot) ++(-0,-0.1) coordinate (bot');
}

\usetikzlibrary{shapes.geometric}
\tikzset{fiber/.style={draw,rectangle,inner sep=2pt,font=\scriptsize}}
\tikzset{pb/.style={draw,regular polygon,regular polygon sides=3,inner sep=0pt,shape border rotate=180,font=\scriptsize}}
\tikzset{pbe/.style={draw,regular polygon,regular polygon sides=3,inner sep=2pt,shape border rotate=180}} 
\tikzset{pbsm/.style={draw,regular polygon,regular polygon sides=3,inner sep=-.5pt,shape border rotate=180,font=\scriptsize}} 
\tikzset{pbflat/.style={draw,regular polygon,regular polygon sides=3,shape border rotate=270,inner sep=1pt,font=\scriptsize}}
\tikzset{pbeflat/.style={draw,regular polygon,regular polygon sides=3,shape border rotate=270,inner sep=2pt}}
\tikzset{pbsmflat/.style={draw,regular polygon,regular polygon sides=3,shape border rotate=270,inner sep=0pt,font=\scriptsize}}
\tikzset{pf/.style={draw,regular polygon,regular polygon sides=3,inner sep=0pt,font=\scriptsize}}
\tikzset{pfe/.style={draw,regular polygon,regular polygon sides=3,inner sep=2pt}} 
\tikzset{pfsm/.style={draw,regular polygon,regular polygon sides=3,inner sep=-.5pt,font=\scriptsize}} 
\tikzset{pfflat/.style={draw,regular polygon,regular polygon sides=3,shape border rotate=90,inner sep=1pt,font=\scriptsize}}
\tikzset{pfeflat/.style={draw,regular polygon,regular polygon sides=3,shape border rotate=90,inner sep=2pt}}
\tikzset{pfsmflat/.style={draw,regular polygon,regular polygon sides=3,shape border rotate=90,inner sep=0pt,font=\scriptsize}}

\tikzset{backwards/.style={z={(-.3,-.7)},x={(-1,0)},y={(0,-1.5)}}}
\tikzset{strings/.style={scale=.75, z={(0,-.4)},x={(-1,0)},y={(0,-2)}}}

\makeatletter

\newif\iftikz@to@relp
\tikz@to@relpfalse
\newif\iftikz@to@relpp
\tikz@to@relppfalse
\tikzoption{abs}[]{\tikz@to@relpfalse\tikz@to@relppfalse}
\tikzoption{rel}[]{\tikz@to@relptrue\tikz@to@relppfalse}
\tikzoption{rrel}[]{\tikz@to@relpfalse\tikz@to@relpptrue}

\tikzstyle{every curve to}=          []
\tikzstyle{curve to}=                [to path=\tikz@to@curve@path]

\tikzoption{bend angle}{\def\tikz@to@bend{#1}}

\tikzoption{bend left}[]{%
  \def\pgf@temp{#1}%
  \ifx\pgf@temp\pgfutil@empty%
  \else%
    \def\tikz@to@bend{#1}%
  \fi%
  \let\tikz@to@out=\tikz@to@bend%
  \c@pgf@counta=180\relax%
  \advance\c@pgf@counta by-\tikz@to@out\relax%
  \edef\tikz@to@in{\the\c@pgf@counta}%
  \tikz@to@switch@on%
  \tikz@to@relativetrue%
}

\tikzoption{bend right}[]{%
  \def\pgf@temp{#1}%
  \ifx\pgf@temp\pgfutil@empty%
  \else%
    \def\tikz@to@bend{#1}%
  \fi%
  \c@pgf@counta=\tikz@to@bend\relax%
  \c@pgf@counta=-\c@pgf@counta\relax%
  \edef\tikz@to@out{\the\c@pgf@counta}%
  \c@pgf@counta=180\relax%
  \advance\c@pgf@counta by-\tikz@to@out\relax%
  \edef\tikz@to@in{\the\c@pgf@counta}%
  \tikz@to@switch@on%
  \tikz@to@relativetrue%
}

\tikzoption{relative}[true]{\csname tikz@to@relative#1\endcsname}
\newif\iftikz@to@relative
\tikz@to@relativefalse

\tikzoption{in}{\def\tikz@to@in{#1}\tikz@to@switch@on}
\tikzoption{out}{\def\tikz@to@out{#1}\tikz@to@switch@on}

\tikzoption{in looseness}{\tikz@to@set@in@looseness{#1}}
\tikzoption{out looseness}{\tikz@to@set@out@looseness{#1}}
\tikzoption{looseness}{\tikz@to@set@in@looseness{#1}\tikz@to@set@out@looseness{#1}}

\tikzoption{in control}{\tikz@to@set@in@control{#1}}
\tikzoption{out control}{\tikz@to@set@out@control{#1}}
\tikzoption{controls}{\tikz@to@parse@controls#1\pgf@stop}

\tikzoption{in min distance}{\tikz@to@set@distances{#1}{}{}{}}
\tikzoption{in max distance}{\tikz@to@set@distances{}{#1}{}{}}
\tikzoption{in distance}{\tikz@to@set@distances{#1}{#1}{}{}}
\tikzoption{out min distance}{\tikz@to@set@distances{}{}{#1}{}}
\tikzoption{out max distance}{\tikz@to@set@distances{}{}{}{#1}}
\tikzoption{out distance}{\tikz@to@set@distances{}{}{#1}{#1}}
\tikzoption{min distance}{\tikz@to@set@distances{#1}{}{#1}{}}
\tikzoption{max distance}{\tikz@to@set@distances{}{#1}{}{#1}}
\tikzoption{distance}{\tikz@to@set@distances{#1}{#1}{#1}{#1}}

\def\tikz@to@set@distances#1#2#3#4{%
  \tikz@to@setifnotempy{#1}{\tikz@to@in@min}{\let\tikz@to@end@compute=\tikz@to@end@compute@looseness}%
  \tikz@to@setifnotempy{#2}{\tikz@to@in@max}{\let\tikz@to@end@compute=\tikz@to@end@compute@looseness}%
  \tikz@to@setifnotempy{#3}{\tikz@to@out@min}{\let\tikz@to@start@compute=\tikz@to@start@compute@looseness}%
  \tikz@to@setifnotempy{#4}{\tikz@to@out@max}{\let\tikz@to@start@compute=\tikz@to@start@compute@looseness}%
  \tikz@to@switch@on%
}

\def\tikz@to@setifnotempy#1#2#3{%
  \def\pgf@temp{#1}%
  \ifx\pgf@temp\pgfutil@empty\else\def#2{#1}#3\fi%
}

\def\tikz@to@set@in@looseness#1{%
  \def\tikz@to@in@looseness{#1}%
  \let\tikz@to@end@compute=\tikz@to@end@compute@looseness%
  \tikz@to@switch@on%
}
\def\tikz@to@set@out@looseness#1{%
  \def\tikz@to@out@looseness{#1}%
  \let\tikz@to@start@compute=\tikz@to@start@compute@looseness%
  \tikz@to@switch@on%
}

\def\tikz@to@parse@controls#1and#2\pgf@stop{\tikz@to@set@in@control{#2}\tikz@to@set@out@control{#1}}

\def\tikz@to@set@in@control#1{%
  \def\tikz@to@in@control{#1}%
  \let\tikz@to@end@compute=\tikz@to@end@compute@control%
  \tikz@to@switch@on%
}
\def\tikz@to@set@out@control#1{%
  \def\tikz@to@out@control{#1}%
  \let\tikz@to@start@compute=\tikz@to@start@compute@control%
  \tikz@to@switch@on%
}

\def\tikz@to@bend{30}

\def\tikz@to@out{45}
\def\tikz@to@in{135}

\def\tikz@to@out@looseness{1}
\def\tikz@to@in@looseness{1}

\def\tikz@to@in@min{0pt}
\def\tikz@to@in@max{10000pt}
\def\tikz@to@out@min{0pt}
\def\tikz@to@out@max{10000pt}

\def\tikz@to@switch@on{\let\tikz@to@path=\tikz@to@curve@path}

\def\tikz@to@curve@path{%
  [every curve to]
  \pgfextra{\iftikz@to@relative\tikz@to@compute@relative\else\tikz@to@compute\fi}
  \tikz@computed@path
  \tikztonodes%
  \pgfextra{\tikz@to@relpfalse\tikz@to@relppfalse}%
}

\def\tikz@to@modify#1#2{%
  \pgfutil@ifundefined{pgf@sh@ns@#1}
  {}%
  {\edef#1{#1.#2}}
}%

\def\tikz@to@compute{%
  \let\tikz@tofrom=\tikztostart%
  \let\tikz@toto=\tikztotarget%
  \tikz@to@modify\tikz@tofrom\tikz@to@out%
  \tikz@to@modify\tikz@toto\tikz@to@in%
  \ifx\tikz@to@start@compute\tikz@to@start@compute@looseness%
    \tikz@to@compute@distance%
  \else%
    \ifx\tikz@from@start@compute\tikz@to@start@compute@looseness%
      \tikz@to@compute@distance%
    \fi%
  \fi%
  \tikz@to@start@compute%
  \tikz@to@end@compute%
  \iftikz@to@relp
    \edef\tikz@computed@path{.. controls \tikz@computed@start and \tikz@computed@end .. +(\tikz@toto)}
  \else
    \iftikz@to@relpp  
      \edef\tikz@computed@path{.. controls \tikz@computed@start and \tikz@computed@end .. ++(\tikz@toto)}
    \else
      \edef\tikz@computed@path{.. controls \tikz@computed@start and \tikz@computed@end .. (\tikz@toto)}
    \fi
  \fi
}

\def\tikz@to@compute@distance{\tikz@scan@one@point\tikz@@to@compute@distance(\tikz@tofrom)}
\def\tikz@@to@compute@distance#1{%
  \def\tikz@first@point{#1}%
  \iftikz@to@relp%
    \tikz@scan@one@point\tikz@@@to@compute@distance([shift={(\tikz@toto)}]\tikz@tofrom)%
  \else%
    \iftikz@to@relpp%
      \tikz@scan@one@point\tikz@@@to@compute@distance([shift={(\tikz@toto)}]\tikz@tofrom)%
    \else%
      \tikz@scan@one@point\tikz@@@to@compute@distance(\tikz@toto)%
    \fi%
  \fi}
\def\tikz@@@to@compute@distance#1{%
  \def\tikz@second@point{#1}%
  \tikz@to@compute@distance@main%
}
\def\tikz@to@compute@distance@main{%
  \pgf@process{\pgfpointdiff{\tikz@first@point}{\tikz@second@point}}%
  \ifdim\pgf@x<0pt\pgf@xa=-\pgf@x\else\pgf@xa=\pgf@x\fi%
  \ifdim\pgf@y<0pt\pgf@ya=-\pgf@y\else\pgf@ya=\pgf@y\fi%
  \pgf@process{\pgfpointnormalised{\pgfqpoint{\pgf@xa}{\pgf@ya}}}%
  \ifdim\pgf@x>\pgf@y%
    \c@pgf@counta=\pgf@x%
    \ifnum\c@pgf@counta=0\relax%
    \else%
      \divide\c@pgf@counta by 255\relax%
      \pgf@xa=16\pgf@xa\relax%
      \divide\pgf@xa by\c@pgf@counta%
      \pgf@xa=16\pgf@xa\relax%
    \fi%
  \else%
    \c@pgf@counta=\pgf@y%
    \ifnum\c@pgf@counta=0\relax%
    \else%
      \divide\c@pgf@counta by 255\relax%
      \pgf@ya=16\pgf@ya\relax%
      \divide\pgf@ya by\c@pgf@counta%
      \pgf@xa=16\pgf@ya\relax%
    \fi%
  \fi%
  \pgf@x=0.3915\pgf@xa%
  \pgf@xa=\tikz@to@out@looseness\pgf@x%
  \pgf@xb=\tikz@to@in@looseness\pgf@x%
  \pgfmathsetlength{\pgf@ya}{\tikz@to@out@min}
  \ifdim\pgf@xa<\pgf@ya%
    \pgf@xa=\pgf@ya%
  \fi%
  \pgfmathsetlength{\pgf@ya}{\tikz@to@out@max}
  \ifdim\pgf@xa>\pgf@ya%
    \pgf@xa=\pgf@ya%
  \fi%
  \pgfmathsetlength{\pgf@ya}{\tikz@to@in@min}
  \ifdim\pgf@xb<\pgf@ya%
    \pgf@xb=\pgf@ya%
  \fi%
  \pgfmathsetlength{\pgf@ya}{\tikz@to@in@max}
  \ifdim\pgf@xb>\pgf@ya%
    \pgf@xb=\pgf@ya%
  \fi%
}

\def\tikz@to@start@compute@looseness{%
  \edef\tikz@computed@start{([shift=(\tikz@to@out:\the\pgf@xa)]\tikz@tofrom)}%
}
\def\tikz@to@end@compute@looseness{%
  \edef\tikz@computed@end{+(\tikz@to@in:\the\pgf@xb)}%
}
\def\tikz@to@start@compute@control{%
  \let\tikz@computed@start=\tikz@to@out@control%
}
\def\tikz@to@end@compute@control{%
  \let\tikz@computed@end=\tikz@to@in@control%
}

\let\tikz@to@start@compute=\tikz@to@start@compute@looseness%
\let\tikz@to@end@compute=\tikz@to@end@compute@looseness%

\def\tikz@to@compute@relative{%
  \tikz@scan@one@point\tikz@@to@compute@relative(\tikztostart)%
}
\def\tikz@@to@compute@relative#1{%
  \def\tikz@tofrom{#1}%
  \tikz@scan@one@point\tikz@@@to@compute@relative(\tikztotarget)%
}
\def\tikz@@@to@compute@relative#1{%
  \def\tikz@toto{#1}%
  \begingroup
    \pgfutil@ifundefined{pgf@sh@ns@\tikztostart}
    {%
      \let\tikz@first@point=\tikz@tofrom%
      \let\tikz@tostart@tikz=\pgfutil@empty
    }%
    {%
      {%
        \tikz@tofrom%
        \pgf@xc=\pgf@x%
        \pgf@yc=\pgf@y%
        {%
          \pgftransformreset%
          \pgftransformshift{\pgfqpoint{\pgf@xc}{\pgf@yc}}%
          \pgftransformrotate{\tikz@to@out}%
          \pgftransformshift{\pgfqpoint{-\pgf@xc}{-\pgf@yc}}%
          \pgf@process{\pgfpointtransformed{\tikz@toto}}%
        }%
        \pgf@xc=\pgf@x%
        \pgf@yc=\pgf@y%
        \pgfpointshapeborder{\tikztostart}{\pgfqpoint{\pgf@xc}{\pgf@yc}}%
        \xdef\tikz@tofrom@smuggle{\noexpand\pgfqpoint{\the\pgf@x}{\the\pgf@y}}
      }%
      \let\tikz@first@point=\tikz@tofrom@smuggle%
      \tikz@first@point%
      \edef\tikz@tostart@tikz{(\the\pgf@x,\the\pgf@y)}%
    }%
    \pgfutil@ifundefined{pgf@sh@ns@\tikztotarget}
    {%
      \let\tikz@second@point=\tikz@toto%
    }%
    {%
      {%
        \tikz@toto%
        \pgf@xc=\pgf@x%
        \pgf@yc=\pgf@y%
        {%
          \pgftransformreset%
          \pgftransformshift{\pgfqpoint{\pgf@xc}{\pgf@yc}}%
          \pgftransformrotate{180}%
          \pgftransformrotate{\tikz@to@in}%
          \pgftransformshift{\pgfqpoint{-\pgf@xc}{-\pgf@yc}}%
          \pgf@process{\pgfpointtransformed{\tikz@tofrom}}%
        }%
        \pgf@xc=\pgf@x%
        \pgf@yc=\pgf@y%
        \pgfpointshapeborder{\tikztotarget}{\pgfqpoint{\pgf@xc}{\pgf@yc}}%
        \xdef\tikz@toto@smuggle{\noexpand\pgfqpoint{\the\pgf@x}{\the\pgf@y}}
      }%
      \let\tikz@second@point=\tikz@toto@smuggle%
    }%
    \tikz@second@point%
    \edef\tikz@totarget@tikz{(\the\pgf@x,\the\pgf@y)}%
    \tikz@to@compute@distance@main%
    \edef\tikz@to@first@distance{\the\pgf@xa}%
    \edef\tikz@to@second@distance{\the\pgf@xb}%
    \pgftransformreset%
    \pgf@process{\tikz@first@point}%
    \pgf@xa=\pgf@x%
    \pgf@ya=\pgf@y%
    \pgf@process{\tikz@second@point}%
    \advance\pgf@x by-\pgf@xa%
    \advance\pgf@y by-\pgf@ya%
    \pgfpointnormalised{}%
    \pgf@xc=\pgf@x%
    \pgf@yc=\pgf@y%
    \pgf@xb=-\pgf@x%
    \pgf@yb=-\pgf@y%
    {%
      \pgftransformshift{\tikz@first@point}%
      \pgftransformcm{\pgf@sys@tonumber\pgf@xc}{\pgf@sys@tonumber\pgf@yc}{\pgf@sys@tonumber\pgf@yb}{\pgf@sys@tonumber\pgf@xc}%
                      {\pgfpointorigin}%
      \pgf@process{\pgfpointtransformed{\pgfpointpolar{\tikz@to@out}{\tikz@to@first@distance}}}%
      \xdef\tikz@computed@start{(\the\pgf@x,\the\pgf@y)}%
    }
    {%
      \pgftransformshift{\tikz@second@point}%
      \pgftransformcm{\pgf@sys@tonumber\pgf@xc}{\pgf@sys@tonumber\pgf@yc}{\pgf@sys@tonumber\pgf@yb}{\pgf@sys@tonumber\pgf@xc}%
                      {\pgfpointorigin}%
      \pgf@process{\pgfpointtransformed{\pgfpointpolar{\tikz@to@in}{\tikz@to@second@distance}}}%
      \xdef\tikz@computed@end{(\the\pgf@x,\the\pgf@y)}%
    }
    \xdef\tikz@computed@path{
      \tikz@tostart@tikz
      .. controls \tikz@computed@start and \tikz@computed@end ..
      \tikz@totarget@tikz}%
  \endgroup
}

\usepackage{cite}

\usepackage{amsmath,amssymb,amsthm}
\usepackage{graphics}
\usepackage{enumerate}
\usepackage{mathtools}
\usepackage{hyperref}
\usepackage{mdwlist}
\usepackage{yhmath}

\let\xto\xrightarrow

\theoremstyle{theorem}

\newtheorem*{thmA}{Theorem A}
\newtheorem*{thmB}{Theorem B}

\newcommand{\dm}{m}
\newcommand{\edm}{{p}}
\newcommand{\hm}{a}
\newcommand{\dn}{n}
\newcommand{\edn}{{q}}
\newcommand{\hn}{b}
\newcommand{\tw}[1]{\widehat{#1}}
\newcommand{\gTop}{{\sf GpTop}}
\newcommand{\Ex}{{\sf Ex}}

\title{Coincidence invariants and 
higher Reidemeister traces }
\date{\today}
\author{Kate Ponto}

\thanks{719 Patterson Office Tower, Department of Mathematics, University of Kentucky, Lexington, KY, kate.ponto@uky.edu\\ 
The author was partially supported by NSF grant DMS-1207670.}

\begin{document}

\maketitle

\begin{abstract}
The Lefschetz number and fixed point index can be thought of as two different descriptions of the same invariant.  
The Lefschetz number is algebraic and defined using homology.  The index is defined more directly from the topology and is a stable homotopy class.
Both the Lefschetz number and index admit generalizations to coincidences and the comparison of these invariants
retains its central role.  In this paper we show that the identification of the Lefschetz number and index using
formal properties of the symmetric monoidal trace extends to coincidence invariants.  This perspective on the coincidence index and Lefschetz number also suggests difficulties 
for generalizations to a coincidence Reidemeister trace.
\end{abstract}

\section*{Introduction}

A {\bf coincidence point} for a pair of maps $f,g\colon M\rightarrow N$ is a point $x$ of $M$ such that $f(x)=g(x)$. 
Coincidence points are a natural generalization of fixed points and there is a corresponding generalization of the Lefschetz fixed point theorem.

\begin{thmA}\label{lefcon} \cite{Lef26}
Suppose $M$ and $N$ are closed, smooth, $\mathbb{Q}$-orientable manifolds of the same dimension 
and $f,g\colon M\rightarrow N$ are continuous maps.
If $f$ and $g$ have no coincidence points then the Lefschetz number of $f$ and $g$ 
\[L(f,g)\coloneqq \sum_i (-1)^i \mathrm{tr}\left( 
\vcenter{
\xymatrix@C=12pt{H_i(M;\bQ)\ar[r]^{f_*}&H_i(N;\bQ)&
H_i(M;\bQ)\\
&H^{\mathrm{dim}(N)-i}(N;\bQ)\ar[r]^-{g^*}\ar[u]^-{-\cap [N]}&H^{\mathrm{dim}(N)-i}(M;\bQ)\ar[u]^{-\cap [M]}
}
}
\right)\] is zero. 
\end{thmA}
The vertical maps above are the Poincar\'{e} duality isomorphism
and they play in essential role in the definition of $L(f,g)$.  The main result of this note is to give a simple proof of 
the following generalization.  
\begin{thmB}\label{lefcong}
Suppose $M$ and $N$ are closed, smooth manifolds and
\[\theta\colon T\nu_{\triangle\subset N\times N}\wedge K\rightarrow L\wedge M_+ \] is  a  stable map for spaces (or spectra) $K$ and $L$. 
If continuous maps $f,g\colon M\rightarrow N$
 have no coincidence points then 
\[\sum_i (-1)^i \mathrm{tr}\left( 
\xymatrix{\tilde{H}_i(M_+\wedge K;\bQ)\ar[r]^-{(f\times g)_*}&\tilde{H}_i(T\nu_{\triangle\subset N\times N}\wedge K;\bQ)\ar[r]^-{\theta_*}&
\tilde{H}_i(L\wedge M_+;\bQ)\\
}
\right)\] is zero. 
\end{thmB}  

While the formulation of this result using the map $\theta$ is nonstandard, Theorem A and  the generalizations in \cite{sav99,sav01} follow from this result.  It
 will allow us to prove a generalization of Theorem A where fundamental classes are replaced by arbitrary homology classes. See \S\ref{lefnumsec}.

The proofs here use duality and trace in symmetric monoidal categories \cite{DP, smc}.
This allows for short, conceptual proofs that are very similar to the corresponding proofs of the Lefschetz fixed 
point theorem \cite{DP} and Reidemeister trace \cite{thesis}.

We finish by considering a generalization of the Lefschetz fixed point theorem for coincidences to Reidemeister traces.  The approach here  does not appear to generalize to 
this case and this 
 failure  is very suggestive.   If there are 
algebraic generalizations of the Reidemeister trace to coincidences similar to the original definition in \cite{Husseini}
we would expect to see an easy comparison of this  invariant  and a  topological invariant using functoriality as for the Lefschetz number.  Instead we 
have a very natural description of the Reidemeister trace for fixed points in terms of the categorical trace 
while the generalization to coincidences suggested by \cite{Crabb10, KW} is fundamentally incompatible with the trace.

\begin{rmk*}
In this paper we focus on closed smooth manifolds.  Many of the results could also be stated 
in terms of compact ENRs (or finite CW complexes) by replacing normal bundles  by mapping cylinders.
\end{rmk*}

\section{Traces for symmetric monoidal categories}

The trace in symmetric monoidal categories is a generalization of
the trace in linear algebra that retains many of the important properties.  In particular, it
satisfies a generalization of invariance of basis and is functorial.  The generalized
trace is a trace for endomorphisms of modules over a commutative ring, 
endomorphisms of chain complexes of modules over a commutative ring, and endomorphisms
of closed smooth manifolds or compact ENRs.  This section is a summary of \cite{DP,LMS,smc}.

Let $\sV$ be a symmetric monoidal category with monoidal product $\otimes$, unit $S$, and symmetry
isomorphism $\gamma$.  

\begin{defn} An object $A$ in $\sV$ is {\bf dualizable} with {\bf dual} $A^\star$ if there are morphisms 
  \[\eta\colon S\rightarrow A\otimes A^\star \text{ and }\epsilon\colon A^\star\otimes A\rightarrow S\] such that the 
  composites
   \[\xymatrix@R=3pt{A\cong S\otimes A\ar[r]^-{\eta\otimes\id}&A\otimes A^\star\otimes A
     \ar[r]^-{\id \otimes \epsilon}&A\otimes S\cong A\\
    A^\star\cong A^\star\otimes S\ar[r]^-{\id\otimes \eta}&A^\star\otimes A\otimes A^\star
    \ar[r]^-{\epsilon \otimes \id}&S\otimes A^\star\cong A^\star}\] 
  are identity maps.
\end{defn}

The map $\eta$ is the \textbf{coevaluation} and $\epsilon$ is the \textbf{evaluation}.  We say that $(A,A^\star)$ is a \textbf{dual pair}.

A module over a commutative ring is dualizable if and only if it is finitely generated and projective.  A chain complex of modules
over a commutative ring is dualizable if and only if it is finitely generated and projective in each degree and only finitely many degrees are nontrivial.

We say a space is {\bf dualizable} if its suspension spectrum is dualizable in the stable homotopy category.  

\begin{prop}\cite[III.4.1, III.5.1]{LMS} \label{topduals}
 If $M$ is a closed smooth manifold 
  then $M_+\coloneqq M\amalg \ast$ is  dualizable.  If $T\nu_M$ is the Thom space of the normal bundle of an embedding of $M$ in $\mathbb{R}^{\edm}$, the dual of $M_+$ is  $\Sigma^{-\edm}T\nu_M$
\end{prop}

If $A$ is dualizable with dual $A^\star$ and $f\colon A\to  A$ is a morphism in $\sV$, the dual of $f$, denoted $f^\star$, is the composite 
\[\xymatrix{A^\star\cong A^\star\otimes S\ar[r]^-{\id\otimes \eta}&A^\star\otimes A\otimes A^\star\ar[rr]^-{\id\otimes f\otimes \id}&& A^\star\otimes A\otimes A^\star\ar[r]^-{\epsilon\otimes \id}&
S\otimes A^\star\cong A^\star}\]

\begin{defn}
  If $A$ is dualizable with dual $A^\star$, $P, Q$ are objects of $\sV$ and $f\colon P\otimes A\rightarrow A\otimes Q$ is an morphism in $\sV$, the {\bf trace} of $f$, $\tr(f)$, is the composite 
  \[\xymatrix@R=8pt{ P\cong P\otimes S\ar[r]^-{\id\otimes\eta}&P\otimes A\otimes A^\star\ar[r]^{f\otimes \id}&A\otimes Q\otimes 
  A^\star\ar[d]^\gamma\\&&A^\star\otimes A\otimes Q\ar[r]^-{\epsilon\otimes \id}&S\otimes Q\cong Q.}\]
\end{defn}

The trace of a linear transformation is the usual linear algebra trace.  The trace of a chain map is the alternating sum of the levelwise traces. If  $f$ is an endomorphism 
of a  closed smooth manifold  and $H_*(-\colon \mathbb{Q})$ is the rational homology functor, the trace of $H_*(f\colon\mathbb{Q})$ is the {\bf Lefschetz number} of $f$.  
The trace of $f$ in the stable homotopy category  is the {\bf fixed point index} of $f$.

It is important to note that the trace is independent of the choice of dual, coevaluation, and evaluation.

\begin{prop}\label{thm:funct-pres-dual}
  Let $F\colon \sV\to \sW$ be a strong symmetric monoidal functor and   $A\in \sV$ be dualizable with dual $A^\star$.  Then $F(A)$ is dualizable with dual $F(A^\star)$.

  For any map $f\colon Q\ten A
  \to A\ten P$, we have
  \[F(\tr(f)) = \tr\left(F(Q)\ten F(A)\to F(Q\ten A)\xto{F(f)} F(A\ten P)\to F(A)\ten F(P)\right).
  \]
\end{prop}
  The rational homology functor is a strong symmetric monoidal functor 
and for each stable map $M\to N$ there is an induced map $H_*(M;\mathbb{Q})\to H_*(N,\mathbb{Q})$.   In particular, if $\pi_0^s(S^0)$ is the 
zeroth stable homotopy group of $S^0$, there is a map $\iota\colon \pi_0^s(S^0)\to \Hom(\mathbb{Q}, \mathbb{Q})$.  The map $\iota$ is injective
and there is an isomorphism $\mathbb{Z}\to \pi_0^s(S^0)$ where the image of $1\in \mathbb{Z}$ under $\iota$ is the identity map of $\mathbb{Q}$.
\begin{cor}\label{lef}\cite{DP}
If $f\colon M\to M$ is an endomorphism of a closed smooth manifold
the image of the fixed point index of $f$ under the injection $\iota$ 
is the Lefschetz number of $f$.
\end{cor}

There are two results about the trace that will be used to compare the invariants defined here with previously defined approaches to coincidence invariants.
Both are easily verified using string diagrams \cite{js91,smc}.
\begin{lem}\label{orstring}
Given dualizable objects $A$ and $B$ and isomorphisms $\psi_A\colon A^\star\to A\otimes P$ and $\psi_B\colon Q\otimes B\to B^\star$ the trace of 
\[\xymatrix{Q\otimes A\ar[r]^{1\otimes f}&Q\otimes B\ar[r]^-{\psi_B}&B^\star\ar[r]^{g^\star}&A^\star\ar[r]^-{\psi_A}&A\otimes P}\]
 is the composite 
\[\xymatrix@R=5pt{Q\simeq Q\otimes S\ar[r]^-{1\otimes \eta_A}&Q\otimes A\otimes A^\star\ar[r]^-{1^2\otimes \psi_A}&Q\otimes A\otimes A\otimes P\ar[r]^-{1\otimes g\otimes f\otimes 1}&\,\\Q\otimes 
B\otimes B\otimes P\ar[r]^-{1\otimes \gamma\otimes 1}
&Q\otimes 
B\otimes B\otimes P\ar[r]^-{\psi_B\otimes 1^2}&B^\star\otimes B\otimes P \ar[r]^-{\epsilon_B\otimes 1}&S\otimes P\simeq P}\]
\end{lem}

\begin{figure}\label{stringcomparison}
  \begin{tabular}{m{28mm}m{2mm}m{36mm}m{2mm}m{32mm}}
\begin{tikzpicture}[scale=.8]
     \node[vert](etap) at (1, 2){$\eta_{A}$};
     \node[vert](f) at (0, 1){$f$};
     \node[vert](exb) at (0, 0){$\psi_B$};
     \node[vert](etaq) at (1.5, -1){$\eta_A$};
     \node[vert](g) at (.75, -2){$g$};
     \node[vert](epp) at (0, -3){$\epsilon_{B}$};
     \node[vert](exa) at (2, -4.5){$\psi_A$};
     \node[vert](eppp) at (1, -5.5){$\epsilon_{A}$};

     \draw (exa) to[out=-90, in =0, looseness=.8]node[ed, swap]{$A$} (eppp);
     \draw (etap) to[out=-180, in =90, looseness=.8]node[ed, swap]{$A$} (f) ;
     \draw (f) to[out=-90, in =90, looseness=.8]node[ed, swap]{$B$} (exb);
      \draw (exb) to[out=-90, in =180, looseness=.7] node[ed, swap]{$B^\star$} (epp);
     \draw (etap) to[out=0, in =90, looseness=.8] node[ed]{$A^\star$}(2,1)  to[out=-90, in =90, looseness=.8]  (2.5,-1)  to[out=-90, in =45, looseness=.8]   (1.5,-3.5) to[out=-135, in =180, looseness=.8]  (eppp);
     \draw (etaq)  to[out=0, in =90, looseness=.7] node[ed, near start, swap]{$A^\star$} (exa);
     \draw (etaq) to[out=-180, in =90, looseness=.8] node[ed, swap]{$A$} (g);
     \draw (g)  to[in=0, out =-90, looseness=.8] node[ed]{$B$}  (epp);
     \draw (-1,2) to[in=145, out =-90, looseness=.8] node[ed]{$Q$}  (exb);
     \draw (exa) to[in=90, out =-45, looseness=.8] node[ed]{$P$}  (2.5,-6);

\end{tikzpicture}
&=&
\begin{tikzpicture}[scale=.8]
     \node[vert](etap) at (.75, 0){$\eta_{A}$};
     \node[vert](f) at (0, -1){$f$};
     \node[vert](exb) at (0, -2){$\psi_B$};
     \node[vert](etaq) at (2.5, 2){$\eta_A$};
     \node[vert](g) at (1.5, 1){$g$};
     \node[vert](epp) at (.75, -3){$\epsilon_{B}$};
     \node[vert](exa) at (3.25, .25){$\psi_A$};
     \node[vert](eppp) at (2.5, -1){$\epsilon_{A}$};

     \draw (exa) to[out=-90, in =0, looseness=.8]node[ed, swap]{$A$} (eppp);
     \draw (etap) to[out=-180, in =90, looseness=.8]node[ed, swap]{$A$} (f) ;
     \draw (f) to[out=-90, in =90, looseness=.8]node[ed]{$B$} (exb);
      \draw (exb) to[out=-90, in =180, looseness=.7] node[ed, swap]{$B^\star$} (epp);
     \draw (etap) to[out=0, in =180, looseness=.8]  (eppp);
     \draw (etaq)  to[out=0, in =90, looseness=.7] node[ed, near start]{$A^\star$} (exa);
     \draw (etaq) to[out=-180, in =90, looseness=.8] node[ed, swap]{$A$} (g);
     \draw (g)  to[in=0, out =-90, looseness=.8] node[ed]{$B$}  (epp);
     \draw (-.75,2) to[in=125, out =-90, looseness=.8] node[ed]{$Q$}  (exb);
     \draw (exa) to[in=90, out =-65, looseness=.8] node[ed,swap]{$P$}  (4,-4);
\end{tikzpicture}

&=&
\begin{tikzpicture}[scale=.8]
     \node[vert](f) at (0, -1){$f$};
     \node[vert](exb) at (0, -2){$\psi_B$};
     \node[vert](etaq) at (.75, 2){$\eta_A$};
     \node[vert](g) at (-.25, 1){$g$};
     \node[vert](epp) at (1, -3){$\epsilon_{B}$};
     \node[vert](exa) at (1.75, 1){$\psi_A$};

     \draw (exa) to[out=-90, in =90, looseness=.8]node[ed, swap]{$A$} (f);
     \draw (f) to[out=-90, in =90, looseness=.8]node[ed, swap]{$B$} (exb);
      \draw (exb) to[out=-90, in =180, looseness=.7] node[ed, swap]{$B^\star$} (epp);
     \draw (etaq)  to[out=0, in =90, looseness=.7] node[ed, near start]{$A^\star$} (exa);
     \draw (etaq) to[out=-180, in =90, looseness=.8] node[ed, swap]{$A$} (g);
     \draw (g)  to[in=125, out =-90, looseness=.8] (1.5,-1) to[in=0, out =-65, looseness=.8] node[ed]{$B$}  (epp);     
     \draw (-1,2) to[in=125, out =-90, looseness=.8] node[ed]{$Q$}  (exb);
     \draw (exa) to[in=90, out =-65, looseness=.8] node[ed]{$P$}  (2.25,-4);
\end{tikzpicture}

\end{tabular}
\caption{Proof of \autoref{orstring}}
\end{figure}

Suppose given dualizable objects $A$ and $B$ in a symmetric monoidal category along with maps $\phi_A\colon A\to A\otimes A$ and 
 $\hat{\beta} \colon Q\otimes B \to B^\star$.  For maps $g\colon A\to B$ and $\alpha\colon S\to A$ define 
$g^{\phi,\alpha}\colon B\to A$ to be the composite 
\[Q\otimes B\xto{\hat{\beta}} B^\star\xto{g^\star}A^\star\simeq A^\star\otimes S
\xto{1\otimes \alpha}A^\star\otimes A\xto{1\otimes \phi_A} A^\star\otimes A\otimes A\xto{\epsilon_A\otimes 1}S\otimes A\simeq A.\]
Let $\beta\colon Q\otimes B\otimes B\to S$ be the composite 
$Q\otimes  B\otimes B\xto{1\otimes \gamma} Q\otimes  B\otimes B \xto{\hat{\beta}\otimes 1}B^\star\otimes B\xto{\epsilon_{B}}S$.

\begin{lem}\label{naivecompare}
The trace of $f\circ g^{\phi,\alpha}$ is the composite 
\[Q\xto{1\otimes \alpha}Q\otimes A\xto{1\otimes \phi_A}Q\otimes A\otimes A\xto{1\otimes g\otimes f}Q\otimes B\otimes B\xto{\beta}S.\]
\end{lem}

\begin{figure}\label{stringcomparison}
  \begin{tabular}{m{32mm}m{2mm}m{38mm}m{2mm}m{28mm}}
\begin{tikzpicture}[scale=.8]
     \node[vert](etap) at (1, 2){$\eta_{B}$};
     \node[vert](a) at (0, 1){$\hat{\beta}$};
     \node[vert](etaq) at (2.25, -1){$\eta_A$};
     \node[vert](g) at (1.5, -2){$g$};
     \node[vert](epp) at (.75, -3){$\epsilon_{B}$};
     \node[vert](beta) at (2, -4){$\alpha$};
     \node[vert](phi) at (2, -5.5){$\phi_A$};
     \node[vert](epq) at (1, -6.5){$\epsilon_A$};
     \node[vert](f) at (2, -7.5){$f$};
     \node[vert](eppp) at (2, -10){$\epsilon_{B}$};
     \draw (beta) --node[ed]{$A$} (phi);
     \draw (phi) --node[ed]{$A$} (f) to [out = -90, in=135, looseness=.8] node[ed,swap, near start]{$B$}(2.5,-9.25) to [out = -45, in=0, looseness=.8]  (eppp);
     \draw (phi) --node[ed, swap]{$A$} (epq);
     \draw (etap) to[out=-180, in =90, looseness=.8]node[ed]{$B$} (a);
      \draw (a) to[out=-90, in =165, looseness=.7] node[ed]{$B^\star$} (epp);
     \draw (etap) to[out=0, in =90, looseness=.8] node[ed]{$B^\star$}(2,1)  to[out=-90, in =90, looseness=.8]  (3.5,-1)  to[out=-90, in =45, looseness=.8]   (2,-9) to[out=-135, in =180, looseness=.8]  (eppp);
     \draw (etaq)  to[out=0, in =45, looseness=.7]node[ed, swap]{$A^\star$} (2.5, -2.75)  to[out=-135, in =45, looseness=.9] (1.5, -3.75) to[out=-135, in =180, looseness=.7]   (epq);
     \draw (etaq) to[out=-180, in =90, looseness=.8] node[ed, swap]{$A$} (g);
     \draw (g)  to[in=0, out =-90, looseness=.8] node[ed]{$B$}  (epp);
    \draw (0,3)  to[in=100, out =-90, looseness=.8] node[ed, swap]{$Q$}  (a);

\end{tikzpicture}
&=&

\begin{tikzpicture}[scale=.8]
     \node[vert](etap) at (3, -3){$\eta_{B}$};
     \node[vert](a) at (1.75, -5){$\hat{\beta}$};
     \node[vert](etaq) at (3, 0){$\eta_A$};
     \node[vert](g) at (2, -2){$g$};
     \node[vert](epp) at (2.5, -6){$\epsilon_{B}$};
     \node[vert](beta) at (5, 2.5){$\alpha$};
     \node[vert](phi) at (5, 1){$\phi_A$};
     \node[vert](epq) at (4.5, -1){$\epsilon_A$};
     \node[vert](f) at (5.5, -2){$f$};
     \node[vert](eppp) at (4.5, -4){$\epsilon_{B}$};

     \draw (beta) --node[ed]{$A$} (phi);
     \draw (phi) --node[ed]{$A$} (f)  to [out = -90, in=0, looseness=.8]  (eppp);
     \draw (phi) --node[ed, swap]{$A$} (epq);
     \draw (etap) to[out=-180, in =90, looseness=.8]node[ed, swap]{$B$} (a);
      \draw (a) to[out=-90, in =180, looseness=.7] node[ed, swap]{$B^\star$} (epp);
     \draw (etap)   to[out=0, in =180, looseness=.8]  (eppp);
     \draw (etaq)   to[out=0, in =180, looseness=.7]node[ed]{$A^\star$}  (epq);
     \draw (etaq) to[out=-180, in =90, looseness=.8] node[ed, swap]{$A$} (g);
     \draw (g)  to[in=15, out =-90, looseness=.8] node[ed]{$B$}  (epp);
    \draw (1,3)  to[in=100, out =-90, looseness=.8] node[ed, swap]{$Q$}  (a);

\end{tikzpicture}

&=&

\begin{tikzpicture}[scale=.8]
     \node[vert](a) at (1, -3){$\hat{\beta}$};
     \node[vert](g) at (1, -1){$g$};
     \node[vert](epp) at (2, -4){$\epsilon_{B}$};
     \node[vert](beta) at (2, 3){$\alpha$};
     \node[vert](phi) at (2, 1){$\phi_A$};
     \node[vert](f) at (3, -1){$f$};

     \draw (beta) --node[ed]{$A$} (phi);
     \draw (phi) --node[ed]{$A$} (f)  to [out = -90, in=90, looseness=.8]  (a);
     \draw (phi) --node[ed, swap]{$A$} (g);
      \draw (a) to[out=-90, in =180, looseness=.7] node[ed, swap]{$B^\star$} (epp);
     \draw (g)  to[in=90, out =-90, looseness=.8] node[ed, near end]{$B$} (3,-3)to[in =0, out=-90]  (epp);
    \draw (.5,4)  to[in=120, out =-90, looseness=.8] node[ed, swap]{$Q$}  (a);

\end{tikzpicture}
\end{tabular}
\caption{Proof of \autoref{naivecompare}}
\end{figure}

\section{Lefschetz Numbers}\label{lefnumsec}
Following \cite{Fuller54,Gon99,GJW06,KW, kos04, kos06_2,kos06, sav99,sav01} we  start  from the observation that 
the coincidence points of maps $f, g\colon M\rightarrow N$ 
are the intersection of the diagonal in $N$ with the image of the product
\[f\times g\colon M\rightarrow N\times N.\]  If we use $\nu_{\triangle \subset N\times N}$ to denote the 
normal bundle of the diagonal in $N\times N$ and $T\nu_{\triangle \subset N\times N}$ to denote 
the Thom space of this bundle, for  coincidence free maps the composite
\[\xymatrix{M\ar[r]^-{f\times g}& N\times N \ar[r]& T\nu_{\triangle \subset N\times N}}\] 
where the second map is the Thom collapse will be homotopic to the constant map at the collapse point.   We denote this composite 
by $f\times g$ since context will make the meaning unambiguous.

To define the invariants described in the introduction and prove comparison results we need some additional structure.  In this paper we encode that structure using  a stable map 
\[\theta\colon  K\wedge T\nu_{\triangle\subset N\times N}\rightarrow M_+\wedge L. \]  This is exactly the structure needed to 
define traces as in the previous section and apply \autoref{thm:funct-pres-dual}.
If $f$ and $g$ have no coincidences the composite 
\[\xymatrix{K\wedge M_+\ar[rr]^-{\id_K\wedge (f\times g)}&& K\wedge T\nu_{\triangle\subset N\times N}\ar[r]^-\theta&  M_+\wedge L}\] will be homotopically trivial.

\begin{defn}
The {\bf coincidence index of $f$ and $g$  relative to $\theta$} is the symmetric monoidal trace of the composite
\[\xymatrix{K\wedge M_+\ar[rr]^-{\id_K\wedge (f\times g)}&& K\wedge T\nu_{\triangle\subset N\times N}\ar[r]^-\theta&  M_+\wedge L}\]  
\end{defn}
%
The homotopy class of the index is clearly trivial if $f$ and $g$ have no coincidences or are homotopic to maps without coincidences.

There is also a corresponding \textbf{Lefschetz number}.  It is  the symmetric monoidal trace  of the composite
\[\xymatrix{ 
H_*(K)\otimes H_*(M_+)\ar[rr]^-{\id\otimes (f\times g)_*}&&
H_*(K)\otimes H_*(T\nu_{\triangle \subset N\times N})\ar[r]^-{\theta_*}&
H_*(L)\otimes H_*(M_+)
}\]
where $H_*(-)$ is rational homology or any other homology theory with a K\"unneth isomorphism.

\autoref{thm:funct-pres-dual} implies the following result.
\begin{thmB}\label{compare}The map induced on homology by the coincidence index of $f$ and $g$ relative to $\theta$ is the 
same as the Lefschetz number of $H_*(f)$ and $H_*(g)$ relative to $H_*(\theta)$.
\end{thmB} 
This theorem is the coincidence generalization of the familiar Lefschetz-Hopf result that compares topologically and algebraically defined fixed point invariants.
As the index is zero for any pair of maps that are homotopic to coincidence free maps, this implies the usual statement that the Lefschetz number 
of coincidence free maps is zero.

To demonstrate the value of this approach we now provide two examples of the required comparison map $\theta$.  These maps 
will allow us to recover Theorem A as well as generalizations in \cite{sav99,sav01}.

\subsection{Orientability}
\label{eg:orient}Let 
$k_*$ be a homology theory and suppose $M$ and $N$ are $k_*$-orientable.
If $k$ is the  spectrum associated to  $k_*$ there are Thom isomorphisms \cite[20.5.8]{MS} 
\[\psi_M\colon k\wedge T\nu_M\cong k\wedge \Sigma^{\edm-\dm} M_+\hspace{.5cm}\text{and}\hspace{.5cm}
\psi_N\colon k\wedge T\nu_{\triangle\subset N\times N}\cong k\wedge \Sigma^\dn N_+\]
where $\nu_M$ is the normal bundle of an embedding of $M$ in $\mathbb{R}^\edm$, $\dm$ is the dimension of $M$, 
and $\dn$ is the dimension of $N$.  

The Thom isomorphism induce the familiar homology isomorphisms \[\tilde{k}_i(T\nu_M)\cong \tilde{k}_i(\Sigma^{\edm-\dm} M_+)\cong \tilde{k}_{i-\edm+\dm}( M_+)\hspace{.3cm}
\tilde{k}_i(T\nu_{\triangle\subset N\times N})\cong \tilde{k}_i(\Sigma^\dn N_+)\cong \tilde{k}_{i-\dn}(N_+)\]
and define a map 
\[\xymatrix@R=10pt{
\theta\colon k\wedge S^{\edm-\dm}\wedge T\nu_{\triangle\subset N\times N}\ar[d]^{1\wedge \psi_N}&& k\wedge S^{\edm-2\dm+\dn}\wedge M_+\\
 k\wedge S^{\edm-\dm+\dn}\wedge N_+\ar[r]^-\pi
&k\wedge S^{\edm-\dm+\dn}\ar[r]^-{1\wedge C}&k\wedge S^{-\dm+\dn}\wedge T\nu_M\ar[u]^{1\wedge \psi_M}
}\]
where $\pi$ is the projection map for $N$ and $C$ is the Thom collapse for an embedding of $M$ in $\mathbb{R}^{\edm}$.

\begin{thm}\label{concomp} If $k_*$ has a K\"unneth isomorphism and $M$ and $N$ are closed smooth $k$-orientable manifolds
 the stable homotopy class of the index with respect to this $\theta$ is the same as  \[
\sum_i (-1)^i \mathrm{tr}\left(\vcenter{\xymatrix{
\tilde{k}_i(M_+)\ar[r]^{f_*}&\tilde{k}_i(N_+)\ar[d]^\cong&\tilde{k}_{\dm-\dn}(S^0)\otimes \tilde{k}_{i}(M_+)\\
&\tilde{k}_{i+\edn-\dn}(T\nu_N)\ar[r]^{(g^\star)_*}&\tilde{k}_{i+\edm-\dn}(T\nu_M)\ar[u]^\cong
}}\right).\]
\end{thm} 
 
Since Poincar\'{e} duality is the composite of Spanier-Whitehead duality and the Thom isomorphism, this sum of traces 
agrees with the sum in Theorem A.  As the index of coincidence free maps is trivial this theorem implies Theorem A.

\begin{proof}
The index of $f$ and $g$ is the symmetric monoidal trace of the composite 
\[\xymatrix{
k\wedge S^\dm\wedge M_+\ar[rr]^-{\id\wedge \id\wedge (f\times g)}&&k\wedge  S^\dm\wedge T\nu_{\triangle\subset N\times N}\ar[r]^-\theta
& k\wedge S^\dn\wedge M_+.
}\]
This agrees with the trace of the composite 
\[\xymatrix{
\tilde{k}_*(S^\dm\wedge M_+)\ar[rr]^-{k_*(f\times g)}&&\tilde{k}_*( T\nu_{\triangle\subset N\times N})\ar[r]^-{k_*(\theta)}& \tilde{k}_*(S^\dn\wedge M_+).
} \] 
The trace is invariant under cyclic permutation and the trace is the identity functor on endomorphisms of the unit object \cite{smc}, so the trace of $k_*(\theta)k_*(f\times g)$ is the same as the composite
\[\xymatrix@C=13pt{
k_*(S^\edm)\ar[d]^{C_*} & k_*(S^{\edm-\dm}\wedge N_+\wedge N_+)\ar[r]& k_*(S^{\edm-\dm}\wedge T\nu_{\triangle\subset N\times N})\ar[d]^{\psi_N}&
\\k_*(T\nu_M)\ar[r]^-{\psi_M}& k_*(S^{\edm-\dm}\wedge M_+)\ar[u]^-{f\wedge g}& k_*(S^{\edm-\dm+\dn}N_+)\ar[r]^{\pi} &k_*(S^{\edm-\dm+\dn})}\]

The Thom isomorphisms are compatible with the evaluation map for $N$ and the coevaluation map for $M$ in the
sense that the diagrams 
\[\xymatrix@C=15pt{S^{\edn-\dn}\wedge N_+\wedge N_+\ar[r]\ar[d]^{\psi_N\wedge \id}&S^{\edn-\dn}\wedge T\nu_{\triangle \subset N\times N}\ar[d]^{\psi_N}&S^\edm\ar[d]_C\ar@/^/[dr]^-{\eta_M}\\
T\nu_N\wedge N_+\ar[r]\ar@/_/[dr]_{\epsilon_N}&S^{\edn}\wedge N_+\ar[d]^\pi&T\nu_M\ar[r]_-{\triangle_T}\ar[d]^{\psi_M}&M_+\wedge T\nu_M\ar[d]^{\id\wedge \psi_M}\\
&S^{\edn}&S^{\edm-\dm}\wedge M_+\ar[r]^-\triangle&S^{\edm-\dm}\wedge M_+\wedge M_+}
\] commute.  ($k$'s have been omitted for readability.)
This allows us to rewrite the composite above as 
\[\xymatrix{k_*(S^\edm)\ar[d]^{\eta_M}& k_*(S^{\edm-\dm}\wedge N_+\wedge N_+)\ar[r]^-{\psi_N}&k_*(S^{\edm-\dm-\edn+\dn}\wedge T\nu_{ N}\wedge N_+)\ar[d]^{\epsilon_N}\\
k_*(M_+\wedge T\nu_M)\ar[r]^-{\psi_M}& k_*(S^{\edm-\dm}\wedge M_+\wedge M_+)\ar[u]^-{f\wedge g}& k_*(S^{\edm-\dm+\dn}) }\]

Then \autoref{orstring} implies this composite agrees with the description given in the statement of the theorem.
\end{proof}

\subsection{(Co)Homology Classes}\label{eg:class} 

This proof suggests a further generalization where we replace the fundamental classes by arbitrary homology and cohomology classes.   Let $k$ be a 
multiplicative cohomology theory with multiplication $\mu$.

\begin{prop}\label{gensav}  For classes $\alpha\in k_\hm(M_+)$ and $\beta\in k^\hn(T\nu_{\triangle\subset N\times N})$ 
the map 
\[ \xymatrix{k_*(S^\hm)\ar[r]^-\alpha& k_*(M_+) \ar[r]^-{(f\times g)_*}&k_*(T\nu_{\triangle \subset N\times N})\ar[r]^-{\beta( -)}&k_*(S^{\hn})}\]
is the trace of the composite 
\[\xymatrix@C=15pt{k_*(S^{\edn-\hn}\wedge N)\ar[rr]^-{\hat{\beta}}&&k_*(T\nu_N)\ar[r]^-{(g^\star)_*}&k_*(S^{\edn-\edm}\wedge T\nu_M)\ar@{.>}[d]^-{-\cap \alpha}\\
&&&
k_*(S^{\edn-\hm}\wedge M_+)\ar[rr]^{f_*}&&k_*(S^{\edn-\hm}\wedge N_+)}\]
\end{prop}
The map $\hat{\beta}$ is as in \autoref{naivecompare}.  We abuse notation and let $-\cap \alpha $ applied to $\delta\in k_c(T\nu_M)$
be the composite along the top and right sides of
 the diagram below. 
\[\xymatrix{S^{\hm+c}\ar[r]^-{\delta\wedge \alpha}\ar[d]^{(\triangle\wedge \id)\circ  \alpha}&T\nu_M\wedge k\wedge M_+\wedge k\ar[r]^-{\id\wedge \triangle}& T\nu_M\wedge k\wedge M_+\wedge M_+\wedge k\ar[d]^\gamma\\
S^c\wedge M_+\wedge M_+\wedge k\ar[r]^-{\id^2\wedge \delta\wedge \id}\ar[d]^{\check{\delta}\wedge \id^2}&M_+\wedge M_+\wedge T\nu_M\wedge k\wedge k\ar[r]^\gamma& T\nu_M\wedge M_+\wedge M_+\wedge k\wedge k
\ar[d]^-{\epsilon\wedge \id\wedge \mu}\\S^\edm\wedge k\wedge M_+\wedge k\ar[r]^-\sim&S^\edm\wedge k\wedge M_+\wedge k\ar[r]^{\id^2\wedge \mu}&S^\edm\wedge M_+\wedge k}\]
We let $\check{\delta}$ be the composite $M_+\wedge S^c\xto{\id\wedge \delta}  M_+\wedge T\nu_M\wedge k\simeq
 T\nu_M\wedge M_+\wedge k\xto{\epsilon \wedge\id} S^{\edm}\wedge k$ and observe that 
 the image of $\delta$ under $-\cap \alpha$ is   the usual cap product of $\alpha$ and   $\check{\delta}$.

\begin{proof}
By
 \autoref{naivecompare},  the trace of the composite 
\[\xymatrix{k_*(S^{\edn-\hn}\wedge N)\ar[r]^{\hat{\beta}}&k_*(T\nu_N)\ar[d]^{(g^\star)_*}\\
k_*(T\nu_M\wedge S^{\edn-\edm-\hm}\wedge M_+)\ar[d]&k_*(S^{\edn-\edm}\wedge T\nu_M)\ar[l]_-{\alpha}\ar@{.>}[d]^-{-\cap \alpha}\\
 k_*(T\nu_M\wedge S^{\edn-\edm-\hm}\wedge M_+\wedge M_+)\ar[r]^-{\epsilon\otimes \id}&k_*(S^{\edn-\hm}\wedge M_+)\ar[r]^{f_*}&k_*(S^{\edn-\hm}\wedge N_+)}\]
is the same as $ \xymatrix{k_*(S^\hm)\ar[r]^-\alpha& k_*(M_+) \ar[r]^-{(f\times g)_*}&k_*(T\nu_{\triangle \subset N\times N})\ar[r]^-{\beta(-)}&k_*(S^{\hn})}$.  
\end{proof}
 Classes $\alpha\in k_\hm(M_+)$ and $\beta\in k^\hn(T\nu_{\triangle\subset N\times N})$ 
are associated to stable maps $\alpha\colon S^\hm\to M_+\wedge k$ and $\beta\colon T\nu_{\triangle\subset N\times N}\wedge S^{-\hn}
\to k$ and using these descriptions the map in \autoref{gensav} can also be written as the composite 
\[\xymatrix{S^\hm\ar[r]^-\alpha& M_+\wedge k\ar[r]^-{f\times g}&T\nu_{\triangle \subset N\times N}\wedge k\ar[r]^-\beta&S^{\hn}\wedge k\wedge k\ar[r]^-{\id\wedge \mu}&S^{\hn}\wedge k.}\]

To compare with \cite{sav99,sav01} we restrict to the case of rational homology.
There is a composite 
\begin{equation}\label{alttrace}\xymatrix@R=10pt{\Hom(H_*(N),H_{\ast+i}(N))\\
\Hom(H^*(N),\mathbb{Q})\otimes H_{\ast+i}(N)\ar[u]_-\sim\ar[r]
& H^*(N)\otimes H_{\ast+i}(N)\ar[r]^-{-\cap -} &H_i(N)}\end{equation}
where the unlabeled  map is the universal coefficients isomorphism. 
For $\alpha\in H_\hm(M_+)$ the \textbf{Lefschetz homomorphism} of $f$ and $g$ evaluated at $\alpha$ \cite[5.1]{sav01} is the image of 
\begin{equation}\label{savlef}H_{*-\edn+\dn}(N)\cong H_{*}(T\nu_N)\xto{g^*}H_{*-\edn+\edm}(T\nu_M)\xto{-\cap \alpha}H_{\ast-\edn+\hm}(M)\xto{f_*} H_{\ast-\edn+\hm}(N)\end{equation}
under \ref{alttrace}.

The \textbf{coincidence index} of $f$ and $g$ \cite[4.2]{sav01} is defined to be the composite 
\[H_*(M_+)\xto{f\times g}H_*(T\nu_{\triangle \subset N\times N})\cong H_{*-\dn}(N_+)\]
where the second map is the Thom isomorphism.
\begin{thm}\cite[2.3]{sav99}\cite{sav01}  The Lefschetz homomorphism at $\alpha$ agrees with the coincidence index at $\alpha$.
 \end{thm}

\begin{proof}
We can determine the resulting element of $H_*(N)$ under either the Lefschetz number or coincidence index by evaluating at each element in $H^*(N)$.

For the Lefschetz homomorphism,  first observe that for each $\iota\in H^{\hm-\dn} (N)$ the image of a homomorphism under 
\begin{equation}\label{alttrace2}
\xymatrix@R=10pt{\Hom(H_{\ast-\edn+\dn}(N),H_{\ast-\edn+\hm}(N))\\ H^{\ast-\edn+\dn}(N)\otimes H_{\ast-\edn+\hm}(N)\ar[u]^\sim\ar[r]^-{-\cap -} &H_{\hm-\dn}(N)\ar[r]^-{\iota(-)} &\mathbb{Q}}\end{equation}
is the same as the image under
\[\xymatrix@R=10pt{\Hom(H_{*-\edn+\dn}(N),H_{\ast-\edn+\hm}(N))\\
 H^{*-\edn+\dn}(N)\otimes H_{\ast-\edn+\hm}(N)\ar[r]^-{\id \otimes (\iota\cap -)}\ar[u]^-\sim& H^{*-\edn+\dn}(N)\otimes H_{*-\edn+\dn}(N) \ar[r]^-{ev} &\mathbb{Q}}\]
Using  universal coefficients we can replace $ev$ with the evaluation map for $H_*(N)$ and then the 
 image of a map $\psi\colon H_{*-\edn+\dn}(N)\to H_{\ast-\edn+\hm}(N)$ under \ref{alttrace2} is the symmetric monoidal trace of 
\[H_{\ast-\edn+\dn}(N)\xto{\psi}  H_{\ast-\edn+\hm}(N)\xto{\iota\cap -} H_{\ast-\edn+\dn}(N).\]
Then the Lefschetz homomorphism evaluated at $\alpha$ and $\iota$ is the trace of 
\[H_{*-\edn+\dn}(N)\cong H_{*}(T\nu_N)\xto{(g^\star)_*}H_{*-\edn+\edm}(T\nu_M)\xto{-\cap (f^*(\iota)\cap \alpha)}H_{\ast-\edn+\dn}(M)\xto{f_*} H_{\ast-\edn+\dn}(N).\]

Applying \autoref{gensav} with
$f^*(\iota)\cap \alpha$ as the homology class for $M_+$
we see that the Lefschetz homomorphism evaluated at $\alpha$ and $\iota$ is the composite
\[ \xymatrix{H_*(S^{\dn-\edm})\ar[r]^-{f^*(\iota)\cap \alpha}& H_*(M_+) \ar[r]^-{(f\times g)_*}&H_*(T\nu_{\triangle \subset N\times N})\ar[r]^-{\beta(-)}&H_*(S^{\dn})}\]
where $\beta$ is  the trivialization of the tangent bundle in rational homology followed by the projection.  In particular, the diagram
\[\xymatrix{
N_+\wedge N_+\wedge H\mathbb{Q}\ar[d]^-{\id\wedge \triangle\wedge \id}\ar[rr]&&T\nu_{\triangle \subset N\times N}\wedge H\mathbb{Q}\ar[d]^{\hat{\beta}\wedge \id}\\
N_+\wedge N_+\wedge N_+\wedge H\mathbb{Q}\ar[r]^{\beta\wedge \id^2}&
S^n\wedge H\mathbb{Q}\wedge N_+\wedge H\mathbb{Q}\ar[r]^-\gamma&S^{\dn}\wedge N_+\wedge H\mathbb{Q}\wedge H\mathbb{Q}}\]
where the unlabeled horizontal map is a collapse and $\hat{\beta}$ is the Thom isomorphism commutes.  
Then \autoref{savcomparison} recovers the form of the coincidence index in \cite{sav01}.

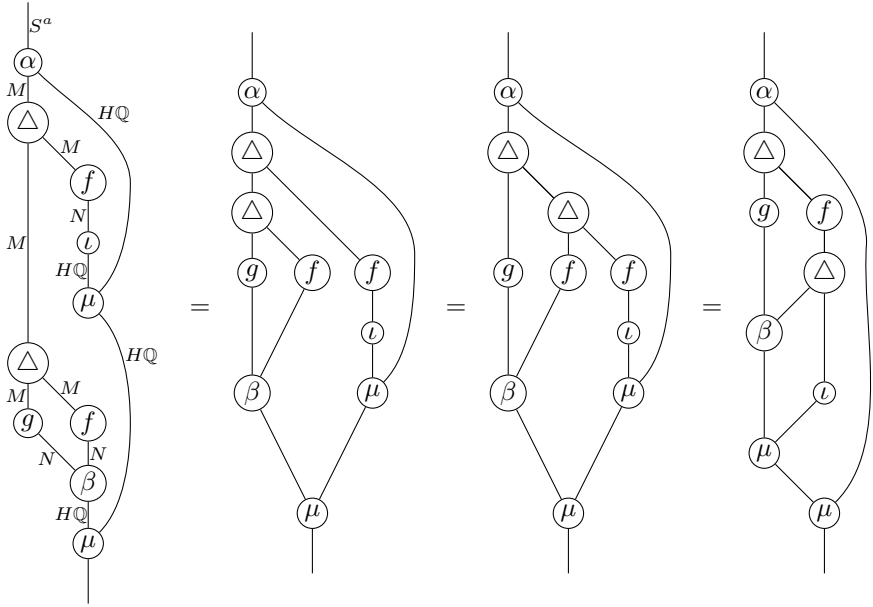
\begin{figure}
  \begin{tabular}{m{21mm}m{2mm}m{25mm}m{2mm}m{25mm}m{2mm}m{18mm}}
\begin{tikzpicture}[scale=.8]
     \node[vert](a) at (0, 0){$\alpha$};
     \node[vert](t1) at (0, -1){$\triangle$};
     \node[vert](f1) at (1, -2){$f$};
     \node[vert](i) at (1, -3){$\iota$};
     \node[vert](mu) at (1, -4){$\mu$};
     \node[vert](t2) at (0, -5){$\triangle$};
     \node[vert](f2) at (1, -6){$f$};
     \node[vert](g) at (0, -6){$g$};
     \node[vert](b) at (1,-7){$\beta$};
      \node[vert](mu2) at (1, -8){$\mu$};

     \draw (0,1) to node[ed]{$S^a$} (a) to node[ed,swap]{$M$} (t1) to node[ed,swap]{$M$} (t2) to node[ed]{$M$} (f2) to node[ed]{$N$} (b);
     \draw (t1) to node[ed]{$M$}  (f1) to node[ed, swap]{$N$} (i) to node[ed, swap]{$H\mathbb{Q}$} (mu);
     \draw (a) to[in=90, out =-45, looseness=.8]  node[ed]{$H\mathbb{Q}$}  (1.7,-2) to[in=45, out =-90, looseness=.8] (mu) ;
    \draw (mu) to[in=90, out =-45, looseness=.8] node[ed]{$H\mathbb{Q}$}  (1.7, -6) to[in=45, out =-90, looseness=.8] (mu2);
      \draw (t2) to  node[ed, swap]{$M$} (g) to  node[ed, swap]{$N$} (b) to node[ed, swap]{$H\mathbb{Q}$}  (mu2) to (1,-9);

\end{tikzpicture}
&=&
\begin{tikzpicture}[scale=.8]
     \node[vert](a) at (0, 0){$\alpha$};
     \node[vert](t1) at (0, -1){$\triangle$};
     \node[vert](f1) at (2, -3){$f$};
     \node[vert](i) at (2, -4){$\iota$};
     \node[vert](mu) at (2, -5){$\mu$};
     \node[vert](t2) at (0, -2){$\triangle$};
     \node[vert](f2) at (1, -3){$f$};
     \node[vert](g) at (0, -3){$g$};
     \node[vert](b) at (0,-5){$\beta$};
      \node[vert](mu2) at (1, -7){$\mu$};

     \draw (0,1) to (a) to (t1) to (t2) to (f2) to (b);
     \draw (t1) to (f1) to (i) to (mu);
     \draw (a) to[in=90, out =-45, looseness=.8]    (2.7,-2.7) to[in=45, out =-90, looseness=.8] (mu) ;
    \draw (mu) to (mu2);
      \draw (t2) to (g) to (b) to (mu2) to (1,-8);
\end{tikzpicture}
&=&
\begin{tikzpicture}[scale=.8]
     \node[vert](a) at (0, 0){$\alpha$};
     \node[vert](t1) at (0, -1){$\triangle$};
     \node[vert](f1) at (2, -3){$f$};
     \node[vert](i) at (2, -4){$\iota$};
     \node[vert](mu) at (2, -5){$\mu$};
     \node[vert](t2) at (1, -2){$\triangle$};
     \node[vert](f2) at (1, -3){$f$};
     \node[vert](g) at (0, -3){$g$};
     \node[vert](b) at (0,-5){$\beta$};
      \node[vert](mu2) at (1, -7){$\mu$};

     \draw (0,1) to (a) to (t1) to (t2) to (f2) to (b);
     \draw (t1) to (t2) to (f1) to (i) to (mu);
     \draw (a) to[in=90, out =-45, looseness=.8]   (2.7,-2.7) to[in=45, out =-90, looseness=.8] (mu) ;
    \draw (mu) to (mu2);
      \draw (t1) to (g) to (b) to (mu2) to (1,-8);
\end{tikzpicture}
&=&
\begin{tikzpicture}[scale=.8]
     \node[vert](a) at (0, 0){$\alpha$};
     \node[vert](t1) at (0, -1){$\triangle$};
     \node[vert](f1) at (1, -2){$f$};
     \node[vert](i) at (1, -5){$\iota$};
     \node[vert](mu) at (0, -6){$\mu$};
     \node[vert](t2) at (1, -3){$\triangle$};
     \node[vert](g) at (0, -2){$g$};
     \node[vert](b) at (0,-4){$\beta$};
      \node[vert](mu2) at (1, -7){$\mu$};

     \draw (0,1) to (a) to (t1) to (f1) to (t2)  to (b);
     \draw (t1) to (f1) to (t2) to (i) to (mu);
     \draw (a) to[in=90, out =-45, looseness=.8]  (1.7,-2.5)  to[in=45, out =-94, looseness=.8]   (mu2);
      \draw (t1) to (g) to (b) to (mu) to  (mu2) to (1,-8);
\end{tikzpicture}
\end{tabular}
\caption{Comparison with the coincidence index in \cite{sav01}.}\label{savcomparison}
\end{figure}
\end{proof}

\subsection{Intersections}
These results are not specific to coincidences but apply to all intersections.  As the generalizations are 
straightforward we only sketch the ideas.

Let $Q$ be a submanifold of a manifold $P$ and $f\colon M\rightarrow P$
be a continuous map.   If the image of $f$ is disjoint from $Q$ the composite
of $f$ with the Thom collapse for the normal bundle of $Q$ in $P$
\[
\xymatrix{M\ar[r]^f&P\ar[r]&T\nu_{Q\subset P}}
\]
is trivial.  
 It is homotopically trivial if $f$ is homotopic to a map $g$ whose image is 
disjoint from $Q$.  In general the converse is not true, see \cite{KW} for a refinement that  
gives a necessary and sufficient condition.  

As in the previous section a stable map $\theta\colon K\wedge T\nu_{Q\subset N}\to L\wedge M_+$ defines both an index 
and Lefschetz number. 

\begin{defn}
The {\bf intersection index of $f$ and $Q$  relative to $\theta$} is the symmetric monoidal trace 
of the composite
\[\xymatrix{K\wedge M_+\ar[r]^-{\id\wedge f}&K\wedge T\nu_{Q\subset N}\ar[r]^-{\theta}&
L\wedge M_+.
}\]  
The {\bf Lefschetz number} is the symmetric monoidal trace of the composite 
\[\xymatrix{H_*(K)\otimes H_*(M_+)\ar[r]^-{\id\otimes f_*}&H_*(K)\otimes H_*(T\nu_{Q\subset N})\ar[r]^-{\theta_*}&
H_*(L)\otimes H_*(M_+)
}\]  where $H_*$ is rational homology.
\end{defn}
With these definitions Theorem B and its proof generalize immediately.

\begin{thm}
The map induced on homology by the coincidence index of $f$ and $Q$ relative to $\theta$ is the 
same as the Lefschetz number of $f$ and $Q$ relative to $H_*(\theta)$.
\end{thm} 

The other examples  for coincidences generalize similarly.

\section{Traces for bicategories}

To extend to the Reidemeister trace we need to replace the trace in a symmetric monoidal 
category with the trace in a bicategory.  This section is a brief summary of the relevant parts of \cite{MS,thesis,shadows}.

\begin{defn} A {\bf bicategory} $\sB$ consists of
  \begin{itemize}\item A collection $\ob\sB$. 
     \item Categories $\sB(A,B)$ for each $A,B\in \mathrm{ob}\sB$. 
    \item Functors \[\odot \colon \sB(A,B)\times \sB(B,C)\rightarrow\sB(A,C)\]
      \[U_A\colon \ast \rightarrow \sB(A,A)\] for $A$, $B$ and $C$ in $\mathrm{ob}\sB$.
  \end{itemize}
  Here $\ast$ denotes the category with one object and one morphism.
  The functors $\odot$ are  required to satisfy unit and
  associativity axioms  up to  natural isomorphisms in $\sB(A,B)$.
\end{defn}

The elements of $\ob\sB$ are called \textbf{0-cells}.  The objects of $\sB(A,B)$
are called \textbf{1-cells}.  The morphisms of $\sB(A,B)$ are called \textbf{2-cells}.

The most suggestive example for the purposes of this paper is the bicategory 
whose objects are rings and for each pair of rings $R$ and  $S$ the associated category is the category of $R$-$S$ bimodules and 
their homomorphisms.      The composition is given by tensor product and a ring regarded as a module over itself is the unit.

\begin{defn}\cite[16.4.1]{MS} A 1-cell $X\in \sB(A,B)$ is {\bf right dualizable} with 
  dual $Y\in \sB(B,A)$ if there are 2-cells \[\xymatrix{\eta\colon U_A\ar[r]&X\odot Y&\epsilon\colon Y\odot X
  \ar[r] &U_B}\] such that  the composites 
  \[\xymatrix@R=5pt{Y\cong Y\odot U_A\ar[r]^-{\id \odot \eta}&
  Y\odot X\odot Y\ar[r]^-{\epsilon\odot \id}&
  U_B\odot Y\cong Y\\
  X\cong U_A\odot X\ar[r]^-{\eta\odot \id}&
  X\odot Y\odot X\ar[r]^-{\id \odot \epsilon}&
  X\odot U_B\cong X}\]
  are identity maps.
\end{defn}

The map $\eta$ is the {\bf coevaluation} and $\epsilon$ is the {\bf evaluation}.  We say $(X,Y)$ is a {\bf dual pair}.

Like the symmetric monoidal trace, the trace of a 2-cell  is defined using a 
composite of the coevaluation and evaluation for a dual pair. 
Unlike that case,  the source of the evaluation and target of the coevaluation are
not isomorphic.   To accommodate this, we need more structure on a bicategory 
before we can define the trace.

\begin{defn}\cite[4.4.1]{thesis}
   A {\bf shadow} for a bicategory $\sB$ is a functor \[\sh{-}\colon \coprod\sB(A,A)\rightarrow \sT\] to a 
   category $\sT$ and  unital and associative natural isomorphisms $\sh{X\odot Y}\cong \sh{Y\odot X}$ for every pair of 
   1-cells $X\in \sB(A,B)$ and $Y\in \sB(B,A)$.
\end{defn}

\begin{defn}\cite[4.5]{thesis}\label{dualdef}
  Let $X$ be a dualizable 1-cell in $\sB$ with dual $Y$ and $f\colon Q\odot X\rightarrow X\odot P$ be a 2-cell in 
  $\sB$.  The {\bf trace} of $f$ is the composite
  \[\xymatrix{{\sh{Q}}\cong \sh{Q\odot U_A}\ar[r]^-{\id \odot \eta}&
  \sh{Q\odot X\odot Y}\ar[d]^-{f\odot \id}\\
  &\sh{X\odot P\odot Y}\ar[r]^\sim&
  \sh{P\odot Y\odot X}\ar[r]^-{\id \odot \epsilon}&\sh{P\odot U_B}\cong\sh{P}.}\] 
\end{defn}

The main result of the next section is a comparison of traces in two different bicategories.  One of these traces is a more familiar description of 
the Reidemeister trace and the other is an alternative description that has some significant technical advantages.  We first describe the relevant bicategories.

\subsection{$\gTop$}
 The 0-cells in the bicategory $\gTop$ are finite groups.   
     A 1-cell $X\colon G\to H$ is a based space with an action of $G\times H$ where the actions of $G$ and $H$ are separately free away from the base point.
    The morphisms from $X\colon G\to H$ to $Y\colon G\to H$ are stable homotopy classes of equivariant maps
    from $X$ to $Y$.

    The bicategorical composition is given by the smash product followed by the quotient by the diagonal action.  
    The unit object associated to a finite group $G$ is $G_+$ regarded as a $G$-$G$ set with a trivial action on the base point. The shadow in $\gTop$ is the quotient by the diagonal action of the group. 

\begin{thm}\cite[3.2.3]{thesis} For a closed smooth manifold  $M$ the universal cover $\tilde{M}_+$
is dualizable as a $\ast\times \pi_1(M)$ space in the bicategory $\gTop$. 

\end{thm}

For a continuous map $f\colon M\to M$ there are induced maps $f_*\colon \pi_1(M)\to \pi_1(M)$ and $\tilde{f}\colon \tilde{M}\to \tilde{M}$ (with some care for base points).
The map $\tilde{f}$ is not $\pi_1(M)$ equivariant but it does induce an equivariant map $\tilde{M}\to \tilde{M}\odot (\pi_1M)_{f}$ 
where $\pi_1(M)_{f_*}$ is the $\pi_1(M)-\pi_1(M)$ set 
$\pi_1(M)$ where the right action is via $f_*$.

\begin{thm}\cite[3.2.3]{thesis}
The trace of $\tilde{f}\colon \tilde{M}\to \tilde{M}\odot (\pi_1M)_{f}$ is the Reidemeister trace.
\end{thm}

\subsection{$\Ex$}

A {\bf parametrized space} over a space 
$B$ is a space $E$ along with maps $\sigma\colon B\rightarrow E$  and $p\colon E\rightarrow B$
such that $p\circ\sigma$ is the identity map of $B$.  A map of parametrized spaces commutes with both $p$ and $\sigma$.
For notation and terminology we will follow \cite{mult} which builds on \cite{MS}.

Parametrized spaces are the 1-cells in a bicategory $\Ex$ defined in \cite{MS}.  The 0-cells are topological spaces and a 
1-cell from $A$ to $B$ is a parametrized space over $A\times B$.  The two cells are fiberwise stable homotopy classes of maps.
  In the examples here  the bicategory composition is given by a fiberwise smash product followed by the pullback along the diagonal map 
and by quotienting the resulting section.   For this bicategory we will follow the notation and conventions of \cite[\S 3]{mult}.

\begin{thm}\cite[18.5.1, 18.6.1]{MS}  If $M$ is a  closed smooth manifold 
  $S^0_M\coloneqq M\amalg M$, regarded as a parametrized space over $\ast\times M$,  is right dualizable.
\label{paradual} \end{thm}

For a map of topological spaces $f\colon X\to Y$, we define spaces 
\[\begin{array}{c}
P(\id,f)\coloneqq \{(\gamma,x)\in Y^I\times X|\gamma(0)=f(x) \}\\P(f,\id)\coloneqq \{(x,\gamma)\in X\times Y^I|\gamma(1)=f(x)\}
\end{array}\]
The first has a map to $Y\times X$ by $(\gamma,x)\mapsto (\gamma(1),x)$ and 
the second has a similar map to $X\times Y$.  These become parametrized spaces with the addition of a disjoint section.  We let $Y_f\coloneqq P(\id,f)\amalg (Y\times X)$ and ${}_fY
\coloneqq P(f,\id)\amalg (X\times Y)$.
\begin{thm}\label{basechangedual}\cite[17.3.1]{MS}
For any map of spaces $f\colon X\to Y$ $({}_fY,Y_f)$ is a dual pair.
\end{thm}

 In the bicategory $\Ex$ the shadow
is given by pulling back along the diagonal map and then quotienting by the resulting section.  In particular,  for an endomorphism $f\colon X\to X$, 
$\sh{X_f}\cong (\Lambda^fX)_+$.  

\section{Reidemeister trace}

Now we consider corresponding generalizations to the Nielsen number and Reidemeister trace to coincidences.  There is a 
coincidence Nielsen number \cite{sch55} but our interest here is  the Reidemeister trace and so we are looking 
for a trace description.  

  A continuous map $f\colon M\to M$ induces a fiberwise map 
$S_M\to S_M\odot M_f$.

\begin{thm}\label{bicatrei}
The bicategorical trace of $S_M\to S_M\odot M_f$ is the Reidemeister trace of $f$.
\end{thm}

There are two ways to approach \autoref{bicatrei}.  We can think of the invariants defined in \cite{KW}
as the definition of the Reidemeister trace and then the identification we require can be found 
in \cite{Crabb10}.  Alternatively, we can use a more classical 
description of the Reidemeister trace in terms of fixed point indices and fixed point classes and apply techniques from \cite{thesis}.  
We give a proof using the second approach here.

\begin{proof} 
  The universal cover $\tilde{M}\to M$ is classified by a map $\phi\colon M\to B\pi_1(M)$ and so 
the $\pi_1(M)$-space $\tilde{M}$ is equivalent to the pullback $M\times_{B\pi_1(M)} E\pi_1(M)$.  Using the notation above
we write this as  \[\tilde{M}_+\cong 
S^0_M\odot {}_\phi B(\pi_1(M))\odot \tw{(E\pi_1(M),\rho)}_+\]
where $   \tw{(E\pi_1(M),\rho)}_+$ is the parametrized space $E\pi_1(M)\amalg B\pi_1(M)\to B\pi_1(M)$ regarded as a space over $B\pi_1(M)\times \ast$.
The space $E\pi_1M$ also has an action of $\pi_1M$ that commutes with the quotient map.

Both $S^0_M$ and ${}_\phi B\pi_1(M)$ are dualizable.  
For $\tw{(E\pi_1(M),\rho)}_+$
we do not have a dual pair in a bicategory, but we do have a map 
\[\triangle_!S^0_{B\pi_1(M)}\to\tw{(E\pi_1(M),\rho)}_+\wedge_{\pi_1(M)} (E\pi_1(M),\rho)_+\] over $B\pi_1(M)\times B\pi_1(M)$ 
and a $\pi_1(M)\times \pi_1(M)$-equivariant map \[(E\pi_1(M),\rho)_+\odot \tw{(E\pi_1(M),\rho)}_+ \to \pi_1(M)_+\]
which make the usual triangle diagrams for a dual pair  commute.
The first map is defined by lifting any path in $B\pi_1(M)$ to $E\pi_1(M)$ and then evaluating at the end points.  Since we 
quotient by $\pi_1(M)$ this will be independent of choices.  For the second map two points in the same fiber 
are taken to the group element that transforms one to the other.

If $\hat{f}\colon B\pi_1(M)\to B\pi_1(M)$ is the map induced by $f$, the commutative diagram
\[\xymatrix{
M\ar[r]^f\ar[d]^\phi&M\ar[d]^\phi\\
B\pi_1(M)\ar[r]^{\hat{f}}& B\pi_1(M)
}\]
defines a map $S_f\odot {}_\phi B\pi_1(M)
\to {}_\phi B\pi_1(M) \odot B\pi_1(M)_{B(\hat{f})}$ \cite[3.3]{mult}.   We can define 
a map \[B\pi_1(M)_{B(\hat{f})}\odot (E\pi_1(M),\rho)_+\to (E\pi_1(M),\rho)_+\wedge_{\pi_1(M)} \pi_1(M)_{f_*}\] by $((\gamma,x), e)\mapsto\tilde{\gamma}(0)$
where $\tilde{\gamma}$ is a lift of $\gamma$ to a path ending at $\hat{f}(e)$.  This is a map over $B\pi_1(M)$ and equivariant with respect to the right action of $\pi_1(M)$.

Using the identification $\tilde{M}_+\cong 
S^0_M\odot {}_\phi B(\pi_1(M))\odot (E\pi_1(M),\rho)_+$
the composite 
\begin{align*}
S^0_M\odot {}_\phi B\pi_1(M)\odot &\tw{(E\pi_1(M),\rho)}_+
\xto{f\odot \id} S^0_M\odot M_f\odot {}_\phi B\pi_1(M)\odot \tw{(E\pi_1(M),\rho)}_+\\
&\xto{\hspace{.6cm}} S^0_M\odot {}_\phi B\pi_1(M) \odot B\pi_1(M)_{B(\hat{f})}\odot \tw{(E\pi_1(M),\rho)}_+\\
&\xto{\hspace{.6cm}} S^0_M\odot {}_\phi B\pi_1(M)\odot \tw{(E\pi_1(M),\rho)}_+\wedge \pi_1(M)_{f_*}.
\end{align*}
is the map $\tilde{f}\colon \tilde{M}\to \tilde{M}\odot (\pi_1M)_{f_*}$ induced by $f$.  
Then the trace of $\tilde{f}\colon \tilde{M}\to \tilde{M}\odot (\pi_1M)_{f_*}$ is the composite \cite[7.5]{shadows} \cite[5.2]{mult}
\[S^0\xto{\tr(f)}\sh{M_f}\xto{} \sh{B\pi_1(M)_{B(\hat{f})}}\xto{}\sh{\pi_1(M)_{f_*}}.\]
The composite of the second and third maps takes a twisted loop in $M$ to 
its associated fixed point class. 
\end{proof}

We now attempt to mimic this description for coincidences. 
If $Q$ is a submanifold of $P$ $S^{\nu_{Q\subset P}}$ is the fiberwise one point compactification 
of the normal bundle of this embedding.  This is a parametrized space over $Q$ where the map is induced by the projection map for the bundle.  The map 
$Q\to S^{\nu_{Q\subset P}}$  is the section at infinity.

Corresponding to the classical Thom collapse there is a fiberwise  homotopy Pontryagin-Thom collapse for $\triangle$ in $N\times N$ 
\[\psi\colon S^0_{N\times N}\rightarrow S^{\nu_{\triangle \subset N\times N}}\odot{} _{i_\triangle}(N\times N)\]\cite[\S 6]{Crabb10} and \cite[II.12]{CJ}.
Composing the fiberwise map induced by $f$ and $g$ $S^0_M\rightarrow S^0_{N\times N}\odot (N\times N)_{f\times g}$ with the homotopy Pontryagin-Thom collapse we have  a map
\begin{equation}\label{secondidx}S^0_M\rightarrow S^{\nu_{\triangle \subset N\times N}}\odot{} _{i_\triangle }(N\times N)\odot (N\times N)_{f\times g}.\end{equation}
Further, this is precisely the invariant that detects intersections.

\begin{thm}\cite[Theorem 3.4]{KW} 
\label{KWiff}
If $\dim(M)+3\leq 2\dim( N)$ the fiberwise stable homotopy class of \ref{secondidx} is trivial if and only if 
there is are maps $f',g'\colon M\rightarrow N$, homotopic to $f$ and $g$, such that $f'$ and $g'$ have no coincidences.
\end{thm}

The natural composite to consider for the coincidence Reidemeister trace is 
\begin{align*}S^0\to S_M^0\odot S^{\nu_M}\to S^{\nu_{\triangle \subset N\times N}}\odot{} _{i_\triangle }(N\times N)\odot (N\times N)_{f\times g}\odot S^{\nu_M}
\\\xto{?} S_M^0\odot {} _{f}N\odot N_{ g} \odot  S^{\nu_M}
\to \sh{\Lambda^{f,g}N}\end{align*}
where the second to last map would need to be a generalization of the Thom isomorphisms for $M$ and $N$.  
To use the approach of the first section we need to rewrite the last two maps using the evaluation for the dual pair 
$(S^0_M, S^{\nu_M})$.  At this point we encounter the major difference between duality in monoidal categories and in bicategories - 
duality in symmetric monoidal categories is symmetric but it is sided in a bicategory.  There is no adjunction that will allow us to introduce the 
the evaluation as we did in the previous section.  This is a major obstruction to defining generalizations of the Reidemeister trace 
like those in \cite{Husseini} for coincidences and suggests that a very different approach may be needed.

\bibliographystyle{plain.bst}
\bibliography{trace_coin}

\end{document}